\documentclass[12pt]{article}
\usepackage{amsfonts,amsmath,amsxtra}
\usepackage{latexsym}
\usepackage{amssymb}
\usepackage{color}
\def\hybrid{\topmargin 0pt      \oddsidemargin 0pt
        \headheight 0pt \headsep 0pt
        \textwidth 16.5cm
        \textheight 23cm
        \marginparwidth 0.0in
        \parskip 5pt plus 1pt   \jot = 1.5ex}
\catcode`\@=11
\def\marginnote#1{}
\newcount\hour
\newcount\minute
\newtoks\amorpm
\hour=\time\divide\hour by60 \minute=\time{\multiply\hour by60
\global\advance\minute by-\hour}
\edef\standardtime{{\ifnum\hour<12 \global\amorpm={am}%
        \else\global\amorpm={pm}\advance\hour by-12 \fi
        \ifnum\hour=0 \hour=12 \fi
      \number\hour:\ifnum\minute<10 0\fi\number\minute\the\amorpm}}
\edef\militarytime{\number\hour:\ifnum\minute<10 0\fi\number\minute}

\def\draftlabel#1{{\@bsphack\if@filesw {\let\thepage\relax
   \xdef\@gtempa{\write\@auxout{\string
      \newlabel{#1}{{\@currentlabel}{\thepage}}}}}\@gtempa
   \if@nobreak \ifvmode\nobreak\fi\fi\fi\@esphack}
        \gdef\@eqnlabel{#1}}
\def\@eqnlabel{}
\def\@vacuum{}
\def\draftmarginnote#1{\marginpar{\raggedright\scriptsize\tt#1}}

\def\draft{\oddsidemargin -0.1truein
        \def\@oddfoot{\sl WLimit.tex \hfil
        \rm\thepage\hfil\sl\today\quad\militarytime}
        \let\@evenfoot\@oddfoot \overfullrule 3pt
        \let\label=\draftlabel
        \let\marginnote=\draftmarginnote
\def\@eqnnum{{\rm (\theequation)}
\rlap{\kern\marginparsep\tt\@eqnlabel}%
\global\let\@eqnlabel\@vacuum}  }

\newfont{\Bbbb}{msbm7 scaled 1\@ptsize00}
\newcommand{\zs}{\raise-1pt\hbox{$\mbox{\Bbbb Z}$}}

\def\titlepage{\@restonecolfalse\if@twocolumn\@restonecoltrue\onecolumn
     \else \newpage \fi \thispagestyle{empty}\c@page\z@
\def\thefootnote{\fnsymbol{footnote}} }
\def\endtitlepage{\if@restonecol\twocolumn \else  \fi
        \def\thefootnote{\arabic{footnote}}
        \setcounter{footnote}{0}}  
\relax
\hybrid
\makeatletter
\newdimen\normalarrayskip            
\newdimen\minarrayskip               
\normalarrayskip\baselineskip \minarrayskip\jot
\newif\ifold             \oldtrue            \def\new{\oldfalse}
\def\arraymode{\ifold\relax\else\displaystyle\fi}
\def\eqnumphantom{\phantom{(\theequation)}} 
\def\@arrayskip{\ifold\baselineskip\z@\lineskip\z@
     \else
     \baselineskip\minarrayskip\lineskip1\baselineskip\fi}
\def\@arrayclassz{\ifcase \@lastchclass \@acolampacol \or
\@ampacol \or \or \or \@addamp \or
   \@acolampacol \or \@firstampfalse \@acol \fi
\edef\@preamble{\@preamble
  \ifcase \@chnum
     \hfil$\relax\arraymode\@sharp$\hfil
     \or $\relax\arraymode\@sharp$\hfil
     \or \hfil$\relax\arraymode\@sharp$\fi}}
\def\@array[#1]#2{\setbox\@arstrutbox=\hbox{\vrule
     height\arraystretch \ht\strutbox
     depth\arraystretch \dp\strutbox
width\z@}\@mkpream{#2}\edef\@preamble{\halign \noexpand\@halignto
\bgroup \tabskip\z@ \@arstrut \@preamble \tabskip\z@ \cr}%
\let\@startpbox\@@startpbox \let\@endpbox\@@endpbox
  \if #1t\vtop \else \if#1b\vbox \else \vcenter \fi\fi
  \bgroup \let\par\relax
  \let\@sharp##\let\protect\relax
  \@arrayskip\@preamble}
\def\eqnarray{\stepcounter{equation}%
              \let\@currentlabel=\theequation
              \global\@eqnswtrue
              \global\@eqcnt\z@
              \tabskip\@centering              
              \let\\=\@eqncr
              $$%
            \halign to \displaywidth  \bgroup
             \eqnumphantom \@eqnsel
      \hskip\@centering                               
    $\displaystyle  \tabskip\z@ {##}$%
    &\global\@eqcnt\@ne \hskip 2\arraycolsep
         $ \displaystyle  \arraymode{##}$\hfil
    &\global\@eqcnt\tw@ \hskip 2\arraycolsep
         $\displaystyle\tabskip\z@{##}$\hfil
         \tabskip\@centering
    &{##}\tabskip\z@\cr}
\makeatother

\def\IC{\mathbb{C}}

\def\IP{\mathbb{P}}

\def\IR{\mathbb{R}}

\def\CC {\mathcal{C}}

\def\CE {\mathcal{E}}
\def\CF {\mathcal{F}}

\def\CH {\mathcal{H}}

\def\CL {\mathcal{L}}

\def\t {{\theta}}

\def\a {{\alpha}}
\def\b {{\beta}}

\def\la{\lambda}

\def\pr {\partial}

\def\zb {\bar{z}}

\def\Arg{{\mathop{\rm Arg}}}

\def\Ker{{\rm Ker}}
\def\Lie{{\rm Lie}}

\def\ov {{\overline}}

\def\<{\langle}
\def\>{\rangle}
\def\ov{\overline}
\newtheorem{te}{Theorem}

\newtheorem{cor}{Corollary}

\newcommand{\proof}{\noindent {\it Proof}. }
\newcommand\bqa{\begin{eqnarray}}
\newcommand\eqa{\end{eqnarray}}
\def\be{\begin{eqnarray}\new\begin{array}{cc}}
\def\ee{\end{array}\end{eqnarray}}
\def\beq{\begin{equation}}
\def\eeq{\end{equation}}
\def\bse{\begin{subequations}}                
\def\ese{\end{subequations}}
\def\bp{\begin{pmatrix}}
\def\ep{\end{pmatrix}}

\newcounter{pac}

\newcounter{pacc}[subsection]

0

\begin{document}

\title{\bf On a matrix element representation of special functions
  associated with toric varieties\footnote{
Contribution to Proceedings of the Nankai Symposium on Mathematical
Dialogue 2021.}}
\author{A.A. Gerasimov, D.R. Lebedev and S.V. Oblezin}
\date{\today}
\maketitle

\renewcommand{\abstractname}{}

\begin{abstract}
\noindent {\bf Abstract}. We develop representation theory approach
to the study of special functions associated with toric varieties. In particular
we show that the corresponding special functions are
 given by matrix elements of certain non-reductive Lie algebras.
\end{abstract}
\vspace{5 mm}

With toric manifolds and their  duals  one can naturally associate generalized
hypergeometric functions introduced by Gelfand-Kapranov-Zelevinski
 \cite{GKZ}. These functions arise in description of quantum cohomology
 of toric manifolds and its mirror dual (see e.g.  \cite{Ba}).
 In this note we consider a
 special instance  of these functions associated with toric varieties
obtained by Hamiltonian reduction of $\IC^N$ (see e.g. \cite{Au}).
 Precisely let us
fix the  standard symplectic structure on $\IC^N$
 \be
  \Omega=\frac{\imath}{2}\sum_{i=1}^N dz_i\wedge d\zb_i
  =\frac{1}{2}\sum_{i=1}^N   d|z_i|^2\wedge d\Arg(z_i),
 \ee
and consider the Hamiltonian action of the compact torus
$T^n=(S^1)^n$ on $\IC^N$. The corresponding toric manifold
$Z=\IC^N//T^n$ is the Hamiltonian reduction under action of
$T^n$. The torus $T^N=(S^1)^N$ and hence its
subtorus $T^n\subset T^N$ act on $Z$ with the following momentum
maps
 \be
  \IC^N\longrightarrow\Lie(T^N)^*,\qquad \IC^N\longrightarrow\Lie(T^n)^*.
 \ee
Define a natural projection
 \be\label{M}
  m\,:\quad{\rm Lie}(T^N)^*\longrightarrow {\rm Lie}(T^n)^*\,,
 \ee
which is dual to the inclusion $T^n\subset T^N$. Let $\{E_i,\,i\in
I\}$ indexed by $I=\{1,\ldots ,N\}$ be a  basis in
$\Lie(T^N)^*$ and let $\{e_\a,\,\a\in
A\}$ indexed by $A=\{1,\ldots ,n\}$ be a basis in $\Lie(T^n)^*$. Then  we have
explicitly
\be
 m(E_j)=\sum_{\a=1}^n e_\a m_j^\a, \qquad j\in I\,.
\ee
Let $\{v_k,\,k\in\ov{A}\}$ indexed by $\ov{A}=\{1,\ldots,N-n\}$ be a
basis in the kernel $\Ker(m)$ of \eqref{M}:
 \be\label{Mrel}
  v_k=\sum\limits_{i=1}^Nv_k^iE_i\in\Lie(T^N)^*,\qquad
  \sum_{i=1}^Nm^\a_iv_k^i\,=\,0,\qquad\a\in A,\quad
  k\in\ov{A}\,.
 \ee

Let $Z_1,\ldots , Z_N$ be cohomology classes dual to the coordinate
hyperplane divisors in $Z$ given by intersects of the coordinate
hyperplanes in $\IC^N$ with $Z$. Then in the quantum cohomology
$QH^*(Z,q)$ the following relations hold (see \cite{Ba}):
 \be\label{relations}
  \prod_{j=1}^N\,Z_j^{m_j^\a}=q_\a,\qquad \a\in A,
 \ee
 where $q_\a=e^{x^\a}$ are quantum deformation parameters.
These relations are deformations of the classical relations 
for $q_{\a}=0$ corresponding to the sets of hyperplanes  with zero
intersections. 

The $T^N$-equivariant  mirror dual to the toric manifold $Z$ is
defined by the data $(Y_q,\CF,\omega_q)$ where $Y_q$ is an affine
space given by the set of solutions of (\ref{relations}), $\CF$ is
given by
 \be\label{potgiv}
  \CF(t,\lambda)=\sum_{i=1}^N (-t_i+\imath
  \lambda_i\ln t_i),\qquad
  \omega_q=\prod_{\a=1}^n \delta(f_\a)\,\frac{dt_1\wedge \cdots \wedge
    dt_N}{t_1\cdots t_N},
 \ee
where
 \be\label{INTER}
  f_\a=e^{-x^\a}\prod_{j=1}^N\,t_j^{m_j^\a}-1,\qquad
  \qquad \a\in A.
 \ee
The special function associated with these data is then given by the
following period  \be\label{wavefunc}
  \Psi_{\la_1,\ldots,\la_N}(q_1,\cdots ,q_n)
  =\int\limits_{\IR^N_+}\!\!\omega_q\,\,\,
  \prod_{i=1}^Nt_i^{\imath \lambda_i}e^{-t_i}\,. 
 \ee

The integral representation \eqref{wavefunc} allows interpretation
as a matrix element of a suitable Lie algebra. Namely, consider the
semi-direct sum $\CL_N$ of abelian Lie algebra $(\mathfrak{gl}_1)^N$ and
the Heisenberg algebra of dimension $2N+1$ defined by the set generators
$\{\CE_i,\CF_i,\CH_i,\,i\in I\}$, central element $\CC$ and relations
 \be \label{LA3}
  [\CE_i,\CH_j]=\delta_{ij}\CE_i, \qquad
  [\CF_i,\CH_j]=-\delta_{ij}\CF_i, \qquad
  [\CE_i,\CF_j]=\CC\delta_{ij}, \qquad i,j\in I.
 \ee
We fix a generic irreducible representation $\pi_{\lambda,c}$ in a
suitable space of functions in $t=(t_1,\cdots, t_N)\in \IR_+^N$
with images of the generators in this representation given by
 \be\label{LA3rep}
  E_i=c\pr_{t_i}, \qquad F_i=t_i, \qquad H_i=t_i\pr_{t_i}+\imath
  \lambda_i, \qquad  C=c,
  \ee
  where $\lambda_i\in \IR,\,\, c\in \IR_+$. In the following for
  convenience we take $c=1$.  
Define left and right  vectors by the following conditions
 \be\label{EQQ1}
  E_i|\psi_R\>=-|\psi_R\>,\qquad i\in I\,; \\
  \sum_{i=1}^N
  v_k^iH_i|\psi_L\>=0,\,\,k\in\ov{A},\qquad
  \prod_{j=1}^NF_j^{m_j^\a}|\psi_L\>=|\psi_L\>,\,\,\a\in A,
 \ee
where $\{v_k,k\in\ov{A}\}$ are given by \eqref{Mrel}. Now  introduce the
following matrix element
 \be\label{WT}
\Psi_\lambda(e^y):=\Psi_{\la_1,\ldots,\la_N}(e^{y^1},\ldots,e^{y^N})
  =\<\psi_L|e^{\sum\limits_{j=1}^{N} y^jH_j}|\psi_R\>,
 \ee
with the non-degenerate pairing $\<\,,\,\>$ given by
 \be\label{ScalProd}
  \<a\,,\,b\>\,:=\,\int\limits_{\IR_+^N}\prod_{i=1}^Nd^{\times}t_i\,\,
  \overline{a(t)}\,b(t),\qquad d^{\times}t_i=\frac{dt_i}{t_i}\,.
 \ee

\begin{te} The matrix element  \eqref{WT} allows the following
integral representation
 \be\label{ToricW}
  \Psi_\lambda(e^y)=\int\limits_{\IR_+^N}
  \prod_{i=1}^N d^\times t_i\,
  \prod_{\a=1}^n \delta\Big(e^{-\sum\limits_{j=1}^N y^jm_j^\a}
  \,\prod_{j=1}^N t_j^{m_j^\a}-1\Big)\, \,\prod_{j=1}^N
  t_j^{\imath\lambda_j}\,e^{-t_j}\,\,,
 \ee
and satisfies the equations:
 \be\label{EQQ}
  \Big\{\prod_{j=1}^N\prod_{k=0}^{m_j^{\a}-1}(\imath
  \lambda_j+k-\pr_{y^j})-e^{\sum\limits_{j=1}^N y^jm_j^\a} \Big\}\,
  \Psi_\lambda(e^y)=0, \qquad \a\in A,
 \ee
 \be\label{EQQ2}
  \sum_{j=1}^N v_k^j\pr_{y^j}\cdot
  \Psi_\lambda(e^y)=0,\qquad k\in\ov{A}\,,
 \ee
with $v_k,\,k\in\ov{A}$ defined in \eqref{Mrel}.
\end{te}
\proof Explicitly for the vectors $\psi_{L,R}$  we have
 \be\label{WV2}
  |\psi_R\>=e^{-\sum\limits_{i=1}^N t_i}, \qquad
  |\psi_L\>=\prod_{i=1}^N t_i^{-\imath\lambda_i}\,\,\prod_{\a=1}^n
  \delta\Big(\prod_{j=1}^Nt_j^{m_j^\a}-1\Big).
 \ee
Substitution of \eqref{WV2} into \eqref{WT}  after a simple
  change of variables  leads to \eqref{ToricW}.
For the second assertion, consider the following polynomials of
generators acting  trivially in representation $\pi_{\lambda,c}$
defined by \eqref{LA3rep}
\be
  \CC_\a=\prod_{j=1}^N (F_jE_j)^{m_j^\a} 
  -\prod_{j=1}^N C^{m_j^{\a}}
  (H_j-\imath \lambda_j)^{m_j^\a},\qquad
  \a\in A.
 \ee
Insertion of these  expressions into the matrix element
\eqref{WT}  then leads to the  first set of equations \eqref{EQQ}.
In addition we have the following set of equations \eqref{EQQ2}
obtained by inserting the linear combinations $\sum\limits_{j=1}^N
v_k^j H_j$ into the matrix element \eqref{WT} and  taking into
  account the defining relation \eqref{EQQ1}. $\Box$

\begin{cor} The matrix element \eqref{WT}   effectively depends on a
  smaller number of variables $x^{\alpha}=\sum_{i=1}^{N} y^{i}
  m_{i}^{\alpha}$, $\alpha\in A$   and allows the following representation 
 \be\label{ToricW0}
  \Psi_{\la_1,\ldots,\la_N}(e^{x_1},\ldots,e^{x_n})
  =\int\limits_{\IR_+^N}\prod_{i=1}^N d^\times t_i\,
  \prod_{\a=1}^n \delta\Big(e^{-x^\a}\,\prod_{j=1}^N
  t_j^{m_j^\a}-1\Big)\,\,\,\prod_{j=1}^N   t_j^{\imath\lambda_j}\,e^{-t_j}\,.
  \ee
This function coincides  with the mirror period
    \eqref{wavefunc} and solves the set of equations for $\a\in A$ 
 \be\label{EQQ4}
 \Big\{\prod_{j=1}^N
\prod_{k=0}^{m_j^{\a}-1}
 (\imath   \lambda_j+k-\sum_{\b=1}^n m_j^\b\pr_{x^\b})-e^{x^\a} \Big\}\,
\Psi_{\la_1,\ldots,\la_N}(e^{x_1},\ldots,e^{x_n})=0.
 \ee
\end{cor}

A simple example of the proposed  construction is provided by mirror
correspondence for complex projective space
  $Z=\IP^{\ell}$. Projective vector space allows the standard 
    Hamiltonian reduction  representation 
  $\IP^{\ell}=\IC^{\ell+1}//S^1$,  where the action is defined as 
\be
S^1\,:\quad\IC^{\ell+1}\longrightarrow\IC^{\ell+1},\qquad
e^{\imath\t}:\,(z_1,\ldots ,z_{\ell+1})\longmapsto
(z_1\, e^{\imath\t},\ldots ,z_{\ell+1}\, e^{\imath\t}).
\ee
Therefore we have $N=\ell+1$, $n=1$ and the relation matrix is given
by $||m_j^\a||=(1,\ldots ,1)$. The resulting representation for the
mirror period \eqref{wavefunc} reads
 \be\label{P}
  \Psi_{\la_1,\ldots,\la_{\ell+1}}(e^x)
  =\int\limits_{\IR_+^N}\prod_{i=1}^{\ell+1} d^\times t_i\,
\,\delta\Big(e^{-x}\prod_{j=1}^{\ell+1}t_j-1\Big)\,\,
    \prod_{j=1}^{\ell+1}
    t_j^{\imath\lambda_j}\,e^{-t_j}.
 \ee
By Theorem 1,  \eqref{P}  coincides with the following matrix element
of the Lie algebra
$\CL_{\ell+1}=(\mathfrak{gl}_1)^{\ell+1}\uplus\CH_{2(\ell+1)}^c$ in
the generic representation \eqref{LA3rep}:
 \be\label{PP}
  \Psi_{\la_1,\ldots,\la_{\ell+1}}(e^x)
  =\<\psi_L|e^{x\sum\limits_{i=1}^{\ell+1}a_iH_i}|\psi_R\>, \qquad \sum_i
  a_i=1,
 \ee
where the left and right vectors are defined by
 \be\label{CONSTR1}
 E_i|\psi_R\>=-|\psi_R\>,\qquad i\in\{1,\ldots,\ell+1\}\,;
 \\ (H_k-H_{\ell+1})|\psi_L\>=0,\quad k\in\{1,\ldots,\ell\},\qquad
 F_1\cdots F_{\ell+1}|\psi_L\>=|\psi_L\>.
 \ee
Particular choice of $\{a_i,\,i\in I\}$ is irrelevant due to
\eqref{CONSTR1}. The defining equations for \eqref{P} are obtained
by using the following operator acting by zero in the representation
\eqref{LA3rep}:
\be
  \CC_{\ell+1}=\prod_{i=1}^{\ell+1}
  F_iE_i-C^{\ell+1}\prod_{i=1}^{\ell+1}(H_i-\imath\lambda_i).
\ee
The  resulting equation
\be
  \Big\{\prod_{i=1}^{\ell+1}(\imath\lambda_i-\pr_x)-e^x\Big\}\,
  \Psi_\lambda(e^x)=0, 
\ee
is in agreement with the proper
equivariant counting of holomorphic maps of $\IP^1$ to $\IP^{\ell}$
(see e.g. \cite{GLO}).\\

  {\it  Acknowledgements:} The research of the second author was
supported by RSF grant 16-11-10075. The work of the third author was
partially supported by the EPSRC grant EP/L000865/1.

\noindent {\small {\bf A.A.G.} {\sl Laboratory for Quantum Field Theory
and Information},\\
\hphantom{xxxx} {\sl Institute for Information
Transmission Problems, RAS, 127994, Moscow, Russia};\\
\hphantom{xxxx} {\sl Interdisciplinary Scientific Center
  J.-V. Poncelet (CNRS UMI 2615),  Moscow, Russia;}\\
\hphantom{xxxx} {\it E-mail address}: {\tt anton.a.gerasimov@gmail.com}}\\
\noindent{\small {\bf D.R.L.}
{\sl Laboratory for Quantum Field Theory
and Information},\\
\hphantom{xxxx}  {\sl Institute for Information
Transmission Problems, RAS, 127994, Moscow, Russia};\\
\hphantom{xxxx} {\sl Moscow Center for Continuous Mathematical
Education,\\
\hphantom{xxxx} 119002,  Bol. Vlasyevsky per. 11, Moscow, Russia};\\
\hphantom{xxxx} {\it E-mail address}: {\tt lebedev.dm@gmail.com}}\\
\noindent{\small {\bf S.V.O.} {\sl
 School of Mathematical Sciences, University of Nottingham\,,\\
\hphantom{xxxx} University Park, NG7\, 2RD, Nottingham, United Kingdom};\\
\hphantom{xxxx} {\sl
 Institute for Theoretical and Experimental Physics,
117259, Moscow, Russia};\\
\hphantom{xxxx} {\it E-mail address}: {\tt oblezin@gmail.com}}
\end{document}